\documentclass[11pt,a4paper,english]{amsart}
\usepackage{babel}
\usepackage[latin1]{inputenc}
\usepackage{amsmath}
\usepackage{amsthm}
\usepackage{amsfonts}
\usepackage{indentfirst}
\usepackage{graphicx}

\usepackage{amssymb}
\usepackage[mathscr]{eucal}

\newtheorem{de}{Definition}
\newtheorem{pro}{Proposition}
\newtheorem{cor}{Corollary}
\newtheorem{teo}{Theorem}
\newtheorem{rem}{Remark}
\newtheorem{lem}{Lemma}
\newtheorem{exa}{Example}

\newcommand{\co}{{\mathcal O}}

\newcommand{\gf}{\mathbb{F}}
\newcommand{\gp}{\mathbb{P}}

\newcommand{\gz}{\mathbb{Z}}

\newcommand{\gr}{\mathbb{R}}
\newcommand{\gc}{\mathbb{C}}

\newcommand{\gn}{\mathbb{N}}
\newcommand{\gG}{\mathbb{G}}

\newcommand{\ck}{\mathcal K}

\newcommand{\cp}{{\mathcal P}}
\newcommand{\cL}{{\mathcal L}}

\newcommand{\findemo}{$\ \ \square$}

\renewcommand{\int}{{\rm int}}

\newcommand{\pic}{{\rm Pic}}
\newcommand{\edim}{{\rm edim\;}}

\title{Curves having one place at infinity and linear systems on rational surfaces}
\author{ F.  Monserrat}
\thanks{Supported by Spain Ministry of Education
 MTM2004-00958, GV05/029  and Bancaixa P1-1A2005-08}
\date{}

\begin{document}

\maketitle

\begin{abstract}
Denoting by ${\mathcal L}_d(m_0,m_1,\ldots,m_r)$ the linear system
of plane curves passing through $r+1$ generic points
$p_0,p_1,\ldots,p_r$ of the projective plane with multiplicity $m_i$
(or larger) at each $p_i$, we prove the Harbourne-Hirschowitz
Conjecture for linear systems ${\mathcal L}_d(m_0,m_1,\ldots,m_r)$
determined by a wide family of systems of multiplicities
$\bold{m}=(m_i)_{i=0}^r$ and arbitrary degree $d$. Moreover, we
provide an algorithm for computing a bound of the regularity of an
arbitrary system $\bold{m}$ and we give its exact value when
$\bold{m}$ is in the above family. To do that, we prove an
$H^1$-vanishing theorem for line bundles on surfaces associated with
some pencils ``at infinity''.

\end{abstract}

\section{Introduction}

This paper deals with the problem of computing the dimension of
linear systems on smooth projective surfaces. The main result
provides, for any arbitrary number of generic points in the
projective plane over the field of complex numbers, $\gp^2$, a
wide family of systems of multiplicities for which the
Harbourne-Hirschowitz Conjecture holds. Moreover, we show an
algorithm for computing upper bounds of the regularity of a system
of multiplicities $\bold{m}=(m_0,m_1,\ldots,m_r)$.

Our proofs of these results are based on Section \ref{infinity},
where we give an $H^1$-vanishing theorem for line bundles on those
surfaces $X$ obtained from $\gp^2$ eliminating (by means of
successive blowing-ups) the indeterminacies of the rational map
$f:\gp^2 \cdots \rightarrow \gp^1$ given by certain pencils of plane
curves. These are the pencils ``at infinity'' associated with
rational projective curves of $\gp^2$ that have one place at
infinity and are smooth in their affine parts. The set formed by the
centers of the blowing-ups used to obtain such a surface $X$ turns
out to be a {\it P-sufficient configuration} (this type of
configurations has been introduced and studied in \cite{g-m-bis},
\cite{g-m} and \cite{g-m-2}). This fact, together with the
simplicity and good properties of the effective semigroup of $X$,
leads up to the above mentioned vanishing theorem. Then,
semicontinuity arguments will allow to deduce our main result.

Fixing $r+1$ points $p_0,p_1,\ldots,p_r$ of $\gp^2$ in generic
position and given $r+1$ non-negative integers $m_0,m_1,\ldots,m_r$,
the linear system ${\mathcal L}_d(m_0,m_1,\ldots,m_r)$ of plane
projective curves of fixed degree $d$ having multiplicity $m_i$ (or
larger) at $p_i$ for each $i$, has an expected dimension (attained
when all the conditions being imposed are independent). Those
systems whose dimension is larger than the expected one are called
{\it special}. The Harbourne-Hirschowitz Conjecture intends to give
a description of all special linear systems. Basically, it asserts
that a linear system is special if and only if it has a multiple
fixed component such that its strict transform on the surface
obtained by blowing-up the points $p_0,p_1\ldots,p_r$ is a
$(-1)$-curve (that is, an integral curve with self-intersection
equal to $-1$ and genus zero). This conjecture goes back to B. Segre
\cite{seg}, being reformulated by several authors (see \cite{harb1},
\cite{gim1}, \cite{hir}, \cite{harb4}, \cite{cil}, \cite{cil2}, and
\cite{cil3} for a survey).

Different approaches have been applied to obtain partial results on
the Harbourne-Hirschowitz Conjecture. It has been proved for
$r+1\leq 9$ points (Castelnuovo was the first to deal with these
cases \cite{castel}, although modern proofs are due to Nagata
\cite{nag}, Gimigliano \cite{gim1} and Harbourne \cite{harb4}).
Arbarello and Cornalba \cite{a-c} treated the homogeneous case with
multiplicity $2$ (that is, $m_0=m_1=\ldots=m_r=2)$ using
infinitesimal deformation theory, and Hirschowitz \cite{hir1} proved
the conjecture for the homogeneous case with multiplicity not
greater than 3, by using a specialization technique (the so-called
{\it Horace method}). This result has been generalized by Ciliberto
and Miranda (\cite{cil} and \cite{cil2}) applying a different
degeneration technique, showing that the Harbourne-Hirschowitz
Conjecture is true for the quasihomogeneous case
$m_1=m_2=\ldots=m_r\leq 3$ and $m_0$ arbitrary, and for the
homogeneous case with multiplicity $m$ up to 12 (the cases $13\leq
m\leq 20$ are treated in \cite{cil4} with the same technique and the
help of a computer program). Using a similar approach, Seibert
\cite{sei} proved the conjecture for the quasihomogeneous case with
$m_1=m_2=\ldots=m_r=4$ and, recently, Laface \cite{laf} has done it
for $m_1=m_2=\ldots=m_r=5$. Other advances have been done by Mignon
\cite{mig2} after proving the conjecture when $m_i\leq 4$ for all
$i$, and Évain \cite{ev}, who proves it for the homogeneous case
when the number of points $r+1$ is a power of 4. Also, using a
refinement of the Horace method (the so-called {\it differential
Horace method}), Alexander and Hirschowitz \cite{a-h-1} obtained a
bound $d_0=d_0(m)$ (only depending on $m$) such that, for any $d\geq
d_0$ and any system of multiplicities $(m_0,m_1,\ldots,m_r)$ with
$m_i\leq m$ for all $i$, the linear system ${\mathcal
L}_d(m_0,m_1,\ldots,m_r)$ is non-special. A result which shows that
the conjecture holds whenever there exist sufficiently many small
multiplicities $m_i\leq 4$, at least one of them being $1$, is the
one recently proved by Bunke and Lossen in \cite{los} by applying
the differential Horace method. More recent advances on the subject
are the papers of S. Yang \cite{yang}, who proves the conjecture
when $m_i\leq 7$ for all $i$, and M. Dumnicki and W. Jarnicki
\cite{dum1}, who do so when $m_i\leq 11$ for all $i$; also, in
\cite{dum2} the conjecture is proved for the homogeneous cases with
multiplicity bounded by 42.

Our contribution to the study of linear systems ${\mathcal
L}_d(m_0,m_1,\ldots,m_r)$ is made in Section \ref{s4}. Using the
iterated blowing-ups (introduced by Kleiman in \cite{kle1} and
\cite{kle2}, and also studied in \cite{harb2}, \cite{roe} and
\cite{fdb2}) and results developed in Section \ref{infinity}, we
deduce a sufficient condition for the non-speciality of a linear
system of that type (Theorem \ref{super}). As a consequence, we
determine, for any arbitrary number of points $r+1\geq 2$, a wide
family of systems of multiplicities $(m_0,m_1,\ldots,m_r)$ for which
the special linear systems of the form ${\mathcal
L}_d(m_0,m_1,\ldots,m_r)$ are completely characterized, proving that
the Harbourne-Hirschowitz Conjecture is true for them. This result
has the particularity of providing, for each arbitrary integer
$r\geq 1$, a large set of systems of multiplicities
$\bold{m}=(m_0,m_1,\ldots,m_{r})$ satisfying the
Harbourne-Hirschowitz Conjecture for which the possible $m_i$ are
unbounded. Moreover, we also determine, for whichever $\bold{m}$ in
the above set, the least degree $d$ such that these multiplicities
impose independent conditions to curves of degree $d$ (that is, the
{\it regularity} of $\bold{m}$).

There are many results giving upper bounds of the regularity of a
system of multiplicities (see \cite{gim}, \cite{cat}, \cite{hir},
\cite{bal}, \cite{b-c}, \cite{xu}, \cite{harbc}, \cite{hhf},
\cite{hr}, \cite{dum1} or \cite{harb3} for a survey). In Section
\ref{s43} we introduce a generalization of the algorithm given in
\cite{roe3}, based on our results in Section \ref{infinity},
providing bounds of the regularity which, in many cases, are
better than the existing ones (as far as the author knows).

Every variety $X$ in this article will be considered over the
field of complex numbers $\gc$. Moreover, $K_X$ will denote a
canonical divisor on $X$.

I thank Javier Fernández de Bobadilla for pointing out to me his
geometric-combinatorial proof of Jung's Theorem on factorization of
automorphisms of the plane.

\section{Preliminaries}\label{pre}

\subsection{Configurations}\label{s1}

In this section we summarize some concepts and notations that will
be used throughout the paper. We start with the definition of
configuration.

An {\it ordered configuration} over $\gp^2$ (a {\it configuration}
in the sequel) will be a finite sequence
${\ck}=(p_0,p_1,\ldots,p_n)$ of closed points such that $p_0$
belongs to $X_0:=\gp^2$ and, inductively, if $i\geq 1$ then $p_i$
belongs to the blowing-up $X_i$ of $X_{i-1}$ at $p_{i-1}$. Among
the points of $\ck$ there is a natural partial ordering: $p_i\leq
p_j$ whenever $p_i=p_j$ or the composition of blowing-ups
$X_j\rightarrow X_i$ maps $p_j$ to $p_i$. We will say that $\ck$
is a {\it chain configuration} when $\leq$ be a total ordering.


Denote by $\pi_{\ck}:Z_{\ck} \rightarrow \gp^2$ the morphism given
by the composition of all the successive blowing-ups centered at
the points of $\ck$. Each blowing-up at $p_i$ gives rise to an
exceptional divisor $E_i$ whose total (resp., strict) transform on
$Z_{\ck}$ will be denoted by $E_i^\ck$ (resp., $\tilde{E}^\ck_i$).
In the same way, for each effective divisor $C$ on $X$, $C^\ck$
(resp., $\tilde{C}^\ck$) will be the total (resp., strict)
transform of $C$ on $Z_{\ck}$. Also, for each divisor $D$ on
$Z_{\ck}$, $[D]$ will denote its class in $\pic(Z_{\ck})$. The
system $\{[L^\ck], [E_0^\ck],[E_1^\ck],\ldots, [E_n^\ck]\}$ is a
$\gz$-basis of $\pic(Z_{\ck})$, $L$ denoting a general line on
$\gp^2$.

A point $p_i\in \ck$ is said to be {\it proximate} to another point
$p_j\in \ck$ (in short, $i\rightarrow j$ or $p_i \rightarrow p_j$)
if either $i=j+1$ and $p_i$ belongs to the exceptional divisor
$E_j$, or $i>j+1$ and $p_i$ belongs to the strict transform on $X_i$
of the exceptional divisor $E_j$. The point $p_i$ is said to be a
{\it free} point if it is proximate to, at most, one point of $\ck$;
otherwise, $p_i$ is said to be a {\it satellite} point. The
proximity relation among the points of $\ck$ is an equivalent datum
to a matrix $\bold{P}_{\ck}=(q_{ij})_{0\leq i,j\leq n}$, called {\it
proximity matrix} of $\ck$, and defined as follows: $q_{ij}=1$ if
$i=j$, $q_{ij}=-1$ if $p_i$ is proximate to $p_j$, and $q_{ij}=0$
otherwise. For each $j=0,1,\ldots,n$, the entries of its $j$th
column are the coefficients of the expression of the divisor
$\tilde{E}_j^{\ck}$ as linear combination of the divisors
$E_0^{\ck}, E_1^{\ck}\ldots, E_n^{\ck}$. The proximity relations can
also be represented by means of a combinatorial object, the {\it
proximity graph}. It will be denoted by $\gG(\ck)$ and it is a
labelled graph whose vertices represent the points of $\ck$ and
whose edges join vertices associated with proximate points. Each
vertex is labelled with the subindex $i$ of its associated point
$p_i$. An edge joining $p_j$ and $p_i$ ($i>j$) is a continuous
straight line whenever $p_i$ is a minimal point of $\ck$ (with
respect to the ordering $\leq $) which is proximate to $p_j$, and it
is a dotted curved line otherwise (the label of an edge is
determined by its property of being continuous-straight or
curved-dotted). For the sake of simplicity, when we will depict a
proximity graph, we will not draw those edges which can be deduced
from others. Notice that the subgraph consisting of the vertices and
the continuous edges has a forest structure whose trees are rooted
on the vertices corresponding to those points in the configuration
which lie in $\gp^2$. A proximity graph will be called {\it
unibranched} if it is associated to a chain configuration. The
proximity graph is an equivalent datum either to the Enriques
diagram or the dual graph of $\ck$.


By a {\it system of multiplicities} we mean a finite sequence
$(m_0,m_1,\ldots,m_n)$ of non-negative integers. A {\it weighted
configuration} (resp., {\it weighted proximity graph}) will be a
pair $(\ck,\bold{m})$ (resp., ($\gG(\ck),\bold{m})$), where
$\ck=(p_i)_{i=0}^n$ is a configuration and
$\bold{m}=(m_0,m_1,\ldots,m_n)$ is a system of multiplicities. It
can be seen as a map that assigns, to each point $p_j$ of $\ck$
(resp., to  the corresponding vertex of $\gG(\ck)$), the
non-negative integer (multiplicity) $m_j$. The {\it excesses} of the
weighted configuration $(\ck,\bold{m})$ are defined to be the
integers $\rho_j(\ck, \bold{m}):=m_j-\sum_{k\rightarrow j} m_k$,
$0\leq j \leq n$. Since the excesses only depend on the proximity
relations among the points of the configuration, we can define the
{\it excesses of a given weighted proximity graph} $(\gG, \bold{m})$
as those associated with every weighted configuration
$(\ck,\bold{m})$ such that $\gG(\ck)=\gG$; they will be denoted by
$\rho_j(\gG,\bold{m})$. Similarly, the {\it proximity matrix
associated with a proximity graph} $\gG$, which will be denoted
$\bold{P}_{\gG}$, can be defined in an obvious way. If $\ck$ is a
configuration of $n+1$ points and $\bold{v}=(v_0,v_1,\ldots,v_t)$ is
a system of multiplicities such that $t<n$, the pair $(\ck,
\bold{v})$ (resp., ($\gG(\ck),\bold{v})$) will also be considered a
weighted configuration (resp., weighted proximity graph),
identifying $\bold{v}$ with the sequence of multiplicities of length
$n+1$ obtained adding  $n-t$ zero components to $\bold{v}$, that is,
$( v_0 ,v_1 , \ldots ,v_t ,0,0, \ldots ,0)$.

\subsection{P-sufficient configurations}\label{psuf}

Consider a configuration $\ck=(p_0,p_1,\ldots,p_n)$ and take the
notations of Section \ref{s1}. Denote by $(b_{ij})_{0\leq i,j\leq
n}$ the entries of the matrix $\bold{P}_{\ck}^{-1}$, whose columns
contain the coefficients of the expressions of each divisor
$E_j^{\ck}$ as linear combinations of the divisors
$\tilde{E}_0^\ck,\tilde{E}_1^\ck,\ldots,\tilde{E}_n^\ck$. For each
integer $i$ such that $0\leq i\leq n$ we consider the divisor on
$Z_{\ck}$ defined by $D_i:=\sum_{j=0}^i b_{ij} E_{j}^{\ck}$, which
has the following property: $D_i\cdot \tilde{E}_j^{\ck}$ equals $-1$
if $i=j$ and $0$ otherwise. We define the square symmetric matrix
$G_{\ck}=(g_{ij})_{0\leq i,j\leq n}$ by
$$g_{ij}=-9D_i\cdot D_j -(K_{Z_{\ck}}\cdot D_i)(K_{Z_{\ck}}\cdot
D_j).$$

Given an element $\bold{x}\in \gr^{n+1}$, we set $\bold{x}> 0$
when all the coordinates of $\bold{x}$ are non-negative and at
least one of them is positive. Recall \cite{minors} that an
($n+1$)-dimensional square symmetric matrix $A$ is called to be
{\it conditionally positive definite} if $\bold{x} A \bold{x}^t>0$
for all vector $\bold{x}\in \gr^{n+1}$ such that $\bold{x}>0$.

\begin{de}
{\rm A configuration $\ck$ is called to be {\it P-sufficient} if
the matrix $G_{\ck}$ is conditionally positive definite. }
\end{de}

This type of configurations has been recently introduced in
\cite{g-m-bis} and \cite{g-m} and, in them, it is proved that the
cone of curves of $Z_{\ck}$ is (finite) polyhedral whenever $\ck$ is
a P-sufficient configuration. Recall that the {\it cone of curves}
of a projective surface $X$, which we will denote by $NE(X)_{\gr}$,
is the convex cone of the real vector space
$\pic(X)\otimes_{\gz}\gr$ spanned by the classes of the effective
divisors on $X$.

When  $\ck$ is a chain configuration, checking whether it is
P-sufficient or not is equivalent to checking a single condition
\cite[Cor. 2]{g-m}: $\ck$ is P-sufficient if and only if
$-9D_n^2-(K_{Z_{\ck}}\cdot D_n)^2>0$.

The following result, whose proof can be found in \cite{g-m-2},
provides a property of the surfaces obtained from P-sufficient
configurations which will be useful in Section \ref{infinity}.

\begin{pro}\label{2}
Let $\ck=(p_0,p_1,\ldots,p_n)$ be a P-sufficient configuration and
$D$ an effective divisor on $Z_{\ck}$ such that $D^2\geq 0$ and
$D\cdot \tilde{E}_i^{\ck}\geq 0$ for all $i=0,1,\ldots,n$. Then,
$K_{Z_{\ck}}\cdot D<0$.
\end{pro}

\subsection{Plane curves having one place at
infinity}\label{plane}

With the exception of Proposition \ref{villa} and Corollary
\ref{pernia}, this section is expository and its aim is to
summarize some facts related to plane curves having one place at
infinity. This type of curves has been extensively studied by
several authors (see, for instance, \cite{abh2}, \cite{moh},
\cite{abh}, \cite{sat}, \cite{pink}, \cite{suzuki} or
\cite{fujimoto}).

\begin{de}\label{definition}
{\rm A projective curve $C\hookrightarrow \gp^2$ (which we will
assume that is not a line) is said {\it to have one place along a
line $H \hookrightarrow \gp^2$} if the intersection $C\cap H$ is a
single point $p$ and $C$ is reduced and has only one analytic
branch at $p$. If $H$ is viewed as the line of infinity in the
compactification of the affine plane to $\gp^2$, we say that {\it
$C$ has one place at infinity.} }
\end{de}

Throughout this section, we fix a projective curve $C$ having one
place at infinity, $p$ being the intersection point of $C$ with
the line of infinity $H$. Consider the infinite sequence of
morphisms
$$\cdots \rightarrow X_{i+1} \rightarrow X_i\rightarrow \cdots
\rightarrow X_1\rightarrow X_0:=\gp^2,$$ where $X_1\rightarrow
X_0=\gp^2$ is the blowing-up of $\gp^2$ at $p_0:=p$ and, for each
$i\geq 1$, $X_{i+1}\rightarrow X_i$ denotes the blowing-up of
$X_i$ at the unique point $p_i$ which lies on the strict transform
of $C$ and on the exceptional divisor $E_{i-1}$ created by the
preceding blowing-up.

\subsubsection{$\delta$-sequences}\label{deltasec}

The unique branch at infinity of $C$ corresponds to a normalized
discrete valuation $v$ of the field of rational functions of $C$
over $\gc$. We define the {\it semigroup at infinity} (resp., {\it
Weierstrass semigroup}) associated with $C$, and we denote it by
$\Gamma_C$ (resp., $S_C$), as the subsemigroup of $\gn$ generated
by all the integers of the form $-v(g)$, where $g$ belongs to the
affine $\gc$-algebra $\co_C(C\setminus \{p\})$ (resp., to the
normalization of $\co_C(C\setminus \{p\})$). Obviously, $\Gamma_C$
is contained in $S_C$ and they are equal if and only if
$C\setminus \{p\}$ is an smooth affine curve. Abhyankar and Moh
proved in \cite{abh2} the existence of a positive integer $s$ and
a sequence of positive generators
$\delta_0,\delta_1,\ldots,\delta_s$ of $\Gamma_C$ such that:
\begin{itemize}
\item[(I)] If $d_i=\gcd (\delta_0,\delta_1,\ldots,\delta_{i-1})$,
for $1\leq i\leq s+1$ and $n_i=d_i/d_{i+1}$, $1\leq i\leq s$, then
$d_{s+1}=1$ and $n_i>1$ for $1\leq i\leq s$.

\item[(II)] For $1\leq i\leq s$, $n_i\delta_i$ belongs to the
semigroup generated by $\delta_0,\delta_1,\ldots, \delta_{i-1}$.

\item[(III)] $\delta_0>\delta_1$ and $\delta_i<\delta_{i-1}
n_{i-1}$ for $i=2,3,\ldots,s$.
\end{itemize}
The sequence $\{\delta_i\}_{i=0}^s$ can be obtained from an
equation of the curve $C$ using approximate roots \cite[Chapter
II, Sections 6,7]{abh2}.  We will refer to it as a {\it
$\delta$-sequence associated with $C$}.

Moreover, it turns out that any sequence
$(\delta_0,\delta_1,\ldots,\delta_s)$ satisfying the above
conditions (I), (II) and (III) is a $\delta$-sequence associated
with some curve having one place at infinity, which can be chosen
of degree $\delta_0$ (see, for instance, \cite{sat} or
\cite{pink}).

Associated with the branch at infinity of a curve having one place
at infinity, there is a sequence of Newton polygons
$P_0,P_1,\ldots,P_{g-1}$ which determines the equisingularity
class of that branch \cite[3.4]{cam}. Assume that each $P_i$ is
the segment which joins the points $(0,e_i)$ and $(m_i,0)$,
$e_i,m_i\in \gn$. These Newton polygons can be explicitly
recovered from a $\delta$-sequence
$(\delta_0,\delta_1,\ldots,\delta_s)$ associated with the curve:

If $\delta_0-\delta_1$ does not divide $\delta_0$, then $s=g$ and
\begin{equation*}
e_0=\delta_0-\delta_1, \;\;\; e_i=d_{i+1}
\end{equation*}
\begin{equation*}
m_0=\delta_0, \;\;\; m_i=n_i\delta_i-\delta_{i+1}
\end{equation*}
for $1\leq i\leq s-1$. Otherwise, $s=g+1$ and
\begin{equation*}
e_0=d_2=\delta_0-\delta_1, \;\;\; e_i=d_{i+2}
\end{equation*}
\begin{equation*}
m_0=\delta_0+n_1\delta_1-\delta_2, \;\;\;
m_i=n_{i+1}\delta_{i+1}-\delta_{i+2}
\end{equation*}
for $1\leq i\leq s-2$.

The above equalities are considered and used in \cite{reg} and the
proximity relations among the infinitely near points
$p_0,p_1,\ldots$ can be easily deduced from the $\delta$-sequence,
as we will describe next (see \cite{cam} for complete details):

Define $s_0=k_0=0$ and let $h_i$, $k_t$ and $s_t$ (with $0\leq
i\leq s_g-1$ and $1\leq t\leq g$) be the positive integers
obtained from the following continued fractions:

$$\frac{m_{j-1}}{e_{j-1}}+k_{j-1}= h_{s_{j-1}}  +
\frac{1}{{h_{s_{j-1}+1}  + _{ \ddots + \frac{1}{{h_{s_j-1}  +
\frac{1}{{k_j }}}}} }},$$ for $j=1,2,\ldots,g$. Also, for each $n\in
\{1,2,\ldots,s_g\}$, define $f(n):=k_t-1$ whenever $n=s_t$ for some
$t\in \{1,2,\ldots,g\}$, and $f(n):=h_n$ otherwise. Then, the
proximity relations are the following: $l\rightarrow l-1$ for each
positive integer $l$, and $l\rightarrow \sum_{i=0}^{n-1} h_i-1$ for
each pair $(n,l)$ such that $1\leq n\leq s_g$ and $\sum_{i=0}^{n-1}
h_i<l\leq \sum_{i=0}^{n-1} h_i+f(n)$.

Thus, a $\delta$-sequence associated with a plane curve $C$ having
one place at infinity determines the equisingularity class of the
branch of $C$ at $p$ and, therefore, the proximity graph of
whichever configuration of the form $(p_0,p_1,\ldots,p_l)$ with
$l\in \gn$ (in particular, that of the minimal embedded resolution
of the branch).

\subsubsection{Curves of Abhyankar-Moh-Suzuki type}\label{ams}

In this paper we are mainly interested in a certain class of
curves having one place at infinity: the so-called curves of
Abhyankar-Moh-Suzuki type, which we define next.

\begin{de}
{\rm A plane curve $C$ having one place at infinity is said to be
of {\it Abhyankar-Moh-Suzuki} type ($AMS$ type for short) if it is
rational and smooth in its affine part, that is, $C\setminus H$ is
isomorphic to $\gc$, $H$ being the line of infinity. }
\end{de}

Let $H$ be the line of infinity in $\gp^2$ and identify $\gc^2$ with
$\gp^2\setminus H$. Recall that, by \cite{abh1}, a curve $C$ is of
AMS type if and only if it is the compactification in $\gp^2$ of the
zero locus of a component of a certain automorphism $\phi:
\gc^2\rightarrow \gc^2$. The embedding of $\gc^2$ in $\gp^2$ allows
to extend $\phi$ to a birational transformation
$\tilde{\phi}:\gp^2\rightarrow \gp^2$. The minimal embedded
resolution of the singularity of $C$ is closely related to the
minimal resolution of the indeterminacy of $\tilde{\phi}$, and the
combinatorics of the last one can be described precisely, as we will
show next. For details see \cite{fdb1} or  \cite{fdb2}.

First, we will define an associative operation $\uparrow$ in the
set of unibranched proximity graphs with two or more vertices (see
Figure 1 for an example).

Let $\gf_1$ and $\gf_2$ be two proximity graphs of this type and
assume that $V_1=\{v_0,v_1,\ldots,v_n\}$ (resp.,
$V_2=\{w_0,w_1,\ldots,w_m\}$) is the set of vertices of $\gf_1$
(resp., $\gf_2$) where, if $\leq$ denotes the ordering induced in
$V_{\gf_1}$ (resp., $V_{\gf_2}$) by the natural ordering among the
points of a configuration whose proximity graph is $\gf_1$ (resp.,
$\gf_2$), it holds that $v_0<v_1\ldots <v_n$ (resp.,
$w_0<w_1\ldots <w_m$). The graph $\gf_1 \uparrow \gf_2$ is the
unibranched proximity graph such that:
\begin{itemize}
\item[- ] its set of vertices is $V_{\gf_1}\cup V_{\gf_2}$;

\item[- ] its set of edges is $A\cup \{e_1,e_2 \}$, where $A$ is
the union of the sets of edges of $\gf_1$ and $\gf_2$ and
$e_1,e_2$ are two new edges such that $e_1$ is a continuous
straight line joining $v_n$ and $w_0$, and $e_2$ is a curved
dotted line joining $v_n$ and $w_1$;

\item[- ] the vertex $v_i$ (resp., $w_i$) is labelled with $i$
(resp., $n+i+1$) for each $i$ such that $0\leq i\leq n$ (resp.,
$0\leq i\leq m$).
\end{itemize}

\begin{figure}[h] \label{figu}
\begin{center}
\setlength{\unitlength}{1mm}
\begin{picture}(55,50)

\qbezier[25](10,5)(5,12.5)(10,20)

\put(10,5){\circle*{1}} \put(10,5){\line(0,1){5}}
\put(10,10){\circle*{1}} \put(10,10){\line(0,1){5}}
\put(10,15){\circle*{1}} \put(10,15){\line(0,1){5}}
\put(10,20){\circle*{1}} \put(10,20){\line(0,1){5}}
\put(10,25){\circle*{1}}

\put(12,5){\scriptsize $0$} \put(12,10){\scriptsize $1$}
\put(12,15){\scriptsize $2$} \put(12,20){\scriptsize $3$}
\put(12,25){\scriptsize $4$}

\qbezier[25](20,10)(15,15)(20,20)

\put(20,5){\circle*{1}} \put(20,5){\line(0,1){5}}
\put(20,10){\circle*{1}} \put(20,10){\line(0,1){5}}
\put(20,15){\circle*{1}} \put(20,15){\line(0,1){5}}
\put(20,20){\circle*{1}}

\put(22,5){\scriptsize $0$} \put(22,10){\scriptsize $1$}
\put(22,15){\scriptsize $2$} \put(22,20){\scriptsize $3$}

\qbezier[25](40,5)(35,12.5)(40,20)

\put(40,5){\circle*{1}} \put(40,5){\line(0,1){5}}
\put(40,10){\circle*{1}} \put(40,10){\line(0,1){5}}
\put(40,15){\circle*{1}} \put(40,15){\line(0,1){5}}
\put(40,20){\circle*{1}} \put(40,20){\line(0,1){5}}
\put(40,25){\circle*{1}}

\qbezier[25](40,35)(35,40)(40,45)

\put(40,30){\circle*{1}} \put(40,30){\line(0,1){5}}
\put(40,35){\circle*{1}} \put(40,35){\line(0,1){5}}
\put(40,40){\circle*{1}} \put(40,40){\line(0,1){5}}
\put(40,45){\circle*{1}}

\put(40,25){\line(0,1){5}}

\qbezier[25](40,25)(35,30)(40,35)

\put(42,5){\scriptsize $0$} \put(42,10){\scriptsize $1$}
\put(42,15){\scriptsize $2$} \put(42,20){\scriptsize $3$}
\put(42,25){\scriptsize $4$}

\put(42,30){\scriptsize $5$} \put(42,35){\scriptsize $6$}
\put(42,40){\scriptsize $7$} \put(42,45){\scriptsize $8$}

\put(9,0){$\gf_1$} \put(19,0){$\gf_2$} \put(35,0){$\gf_1 \uparrow
\gf_2$}

\end{picture}
\end{center}
\caption{Operation $\uparrow$}
\end{figure}
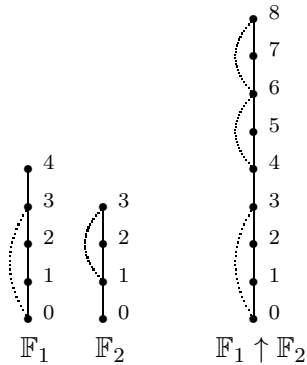

For each integer $n\geq 2$ consider a chain configuration
$(p_0,p_1,\ldots,p_{2n-2})$ such that $p_i$ is proximate to $p_0$
for all $i$ such that $1\leq i\leq n-1$ and the remaining points
$p_{n}, p_{n+1}\ldots p_{2n-2}$ are free. Define $\gG(n)$ to be
the proximity graph of this configuration (see Figure 2).

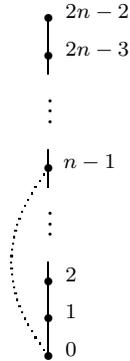
\begin{figure}[h] \label{figu}
\begin{center}
\setlength{\unitlength}{1mm}
\begin{picture}(25,50)

\qbezier[30](15,0)(5,12.5)(15,25)

\put(15,0){\circle*{1}} \put(15,0){\line(0,1){5}}
\put(15,5){\circle*{1}} \put(15,5){\line(0,1){5}}
\put(15,10){\circle*{1}} \put(15,10){\line(0,1){2.5}}

\put(15,25){\circle*{1}} \put(15,25){\line(0,1){2.5}}
\put(15,25){\line(0,-1){2.5}}

 \put(15,40){\line(0,-1){2.5}}

\put(15,40){\circle*{1}} \put(15,40){\line(0,1){5}}
\put(15,45){\circle*{1}}

\put(14.7,16){$\vdots$}


\put(14.7,31){$\vdots$}

\put(16,0){ \scriptsize 0}\put(16,5){ \scriptsize 1}\put(16,10){
\scriptsize 2}\put(16,40){ \scriptsize $2n-3$}\put(16,45){
\scriptsize $2n-2$} \put(17,25){\scriptsize $n-1$}



\end{picture}
\end{center}
\caption{Proximity graph $\gG(n)$}
\end{figure}

For each finite ordered sequence $(n_1,n_2,\ldots,n_r)$ of
integers such that $r\geq 1$ and $n_i\geq 2$ for all $i$, we
define a proximity graph, depending only on that sequence, by
using the above considered associative operation:
$$\gG(n_1,n_2,\ldots,n_r):=\gG(n_1) \uparrow \gG(n_2)
\uparrow \cdots \uparrow \gG(n_r).$$ Also, we denote by
$\gG(n_1,n_2,\ldots,n_r)^{-}$ (resp.,
$\gG(n_1,n_2,\ldots,n_r)^{+}$) the proximity graph obtained from
$\gG(n_1,n_2,\ldots,n_r)$ by deleting (resp., adding) the last
$n_r-1$ vertices and the edges which are adjacent to them (resp.,
a new vertex with label $2\sum_{i=1}^r n_i -r$ and a new edge
joining it with the vertex with label $2\sum_{i=1}^r n_i-r-1$).

Now consider, as above, a curve  $C$ of AMS type and an affine
automorphism $\phi: \gc^2\rightarrow \gc^2$ such that $C$ is the
zero locus of a component of it. Let $\pi:X\rightarrow \gp^2$ be the
minimal resolution of the indeterminacy of
$\tilde{\phi}:\gp^2\rightarrow \gp^2$ and let ${\mathcal K}$ be the
configuration of centers of the blowing-ups involved in $\pi$. Then,
there exists a sequence of integers $(n_1,n_2,\ldots,n_r)$ (with
$n_i \geq 2$ for all $i$) such that $\gG({\mathcal
K})=\gG(n_1,n_2,\ldots,n_r)$. Moreover, the strict transform on $X$
of the line of infinity, $\tilde{H}^{\ck}$, is a $(-1)$-curve, that
is, a smooth rational curve with self-intersection $-1$.

If ${\mathcal C}$ is the configuration such that $\pi_{\mathcal
C}:X_{\mathcal C}\rightarrow \gp^2$ induces the minimal embedded
resolution of the singularity of $C$ at infinity, there are two
possibilities: either $\pi_{\mathcal C}$ is the composition of all
the blowing-ups of $\pi$ except the last $n_r-1$ of them (in this
case, $\gG({\mathcal C})=\gG(n_1,n_2,\ldots,n_r)^-$), or it is the
composition of the first $2\sum_{i=1}^{r-2}n_i+n_{r-1}-r+2$
blowing-ups of $\pi$ (in this case, $\gG({\mathcal
C})=\gG(n_1,n_2,\ldots,n_{r-1})^-$). The following proposition shows
that, for each proximity graph of the form
$\gG(n_1,n_2,\ldots,n_r)^-$, there exists a curve of AMS type such
that the proximity graph associated to its minimal embedded
resolution is this one. Then, the proximity graphs associated to
minimal resolutions of curves of AMS type are exactly those of the
form $\gG(n_1,n_2,\ldots,n_r)^-$.

\begin{pro}\label{villa}
Let $(n_1,n_2,\ldots,n_r)$ be an ordered sequence of integers such
that $n_i\geq 2$ for all $i=1,2,\ldots,r$. Then, there exists a
curve $C$ of AMS type such that its degree is $n_1n_2\cdots n_r$
and the proximity graph associated with its minimal embedded
resolution is $\gG(n_1,n_2,\ldots,n_r)^-$.
\end{pro}
\noindent {\it Proof}. Define the integers $\delta_k=n_{k+1}
n_{k+2}\cdots n_r$ for $k=0,1,\ldots,r-1$ and $\delta_r=1$. It is
obvious that the sequence $(\delta_0,\delta_1,\ldots,\delta_r)$
satisfies the conditions (I), (II) and (III) which characterize
the $\delta$-sequences and, therefore, there exists a curve $C$ of
degree $\delta_0=n_1 n_2 \cdots n_r$ having one place at infinity
with associated $\delta$-sequence
$(\delta_0,\delta_1,\ldots,\delta_r)$. From this sequence one can
compute, using the formulae given in Section \ref{deltasec}, the
proximity relations among the points of the configuration which
provides the minimal embedded resolution of the singularity of $C$
and check that the proximity graph associated with this
configuration is $\gG(n_1,n_2,\ldots,n_r)^-$. Finally, since
$\delta_r=1$, the Weierstrass semigroup of $C$ and its semigroup
at infinity are both equal to $\gn$ and, therefore, $C$ is
rational and smooth in its affine part.\findemo\\

A direct consequence of the above proposition and the genus
formula is the following

\begin{cor}\label{pernia}
Let $C$ be a curve of AMS type and $n_1,n_2,\ldots,n_r\geq 2$
integers such that the proximity graph associated with its minimal
embedded resolution is $\gG(n_1,n_2,\ldots,n_r)^-$. Then, the
degree of $C$ is $n_1 n_2\cdots n_r$.
\end{cor}

\section{Surfaces associated with pencils ``at
infinity''}\label{infinity}

Let $C$ be a projective curve of $\gp^2$ having one place at
infinity and consider the notations of Section \ref{pre}. Take
projective coordinates $(X:Y:Z)$ on $\gp^2$ such that $Z=0$ be the
equation of the line of infinity $H$ and let $F(X,Y,Z)$ be an
homogeneous polynomial in $k[X,Y,Z]$ such that $F(X,Y,Z)=0$ is an
equation of $C$. The {\it pencil ``at infinity'' associated with
$C$}, which we will denote by $\cp(C)$, will be the linear
subspace of $H^0(\gp^2, \co_{\gp^2}(d))$ spanned by $F$ and $Z^d$,
$d$ being the degree of $F$. Let $n$ be the smallest integer such
that the composition of morphisms $X_{n+1}\rightarrow
X_n\rightarrow \cdots \rightarrow X_0=\gp^2$ eliminates the
indeterminacies of the rational map $\gp^2 \cdots \rightarrow
\gp^1$ defined by $\cp(C)$. We will denote by $\ck_C$ the chain
configuration $(p_0,p_1,\ldots,p_n)$ and by $X_C$ the surface
$Z_{\ck_C}=X_{n+1}$. It turns out that all the curves of $\cp(C)$,
except the non-reduced one with equation $Z^d=0$, are integral
curves having one place at infinity, and the above morphism
$X_{C}\rightarrow \gp^2$ induces a simultaneous embedded
resolution of all of them (see \cite{moh}).

The objective of this section is to give a vanishing theorem for
line bundles on the surface $X_C$, when $C$ is a curve of AMS
type. This result will allow to determine the dimension of
whichever complete linear system on $X_C$.

The following proposition provides two characterizations of the
curves of AMS type depending on the associated configuration
$\ck_C$.

\begin{pro}\label{rat}
Let $C$ be a curve having one place at infinity. Then, the
following conditions are equivalent:
\begin{itemize}
\item[(a)] The configuration $\ck_C$ is P-sufficient.

\item[(b)] $K_{X_C}\cdot \tilde{C}^{\ck_C}<0$.

\item[(c)] $C$ is a curve of AMS type.
\end{itemize}
\end{pro}

\noindent {\it Proof}. First, observe that the class of
$\tilde{C}^{\ck_C}$ in $\pic(X_C)$ coincides with
$[dL^{\ck_C}-D_n]$, $d$ being the degree of $C$ and $D_n$ the
divisor defined in Section \ref{psuf} (associated with the
configuration $\ck_C$). The reason is that $D_n=\sum_{i=0}^n
u_iE_i^{\ck_C}$, $u_i$ being the multiplicity of the strict
transform of $C$ at the point $p_i$, $0\leq i\leq n$, because
$D_n\cdot \tilde{E}_i^{\ck_C}$ equals $-1$ if $i=n$ and $0$
otherwise.

The equivalence between $(a)$ and $(b)$ is consequence of the
equalities
$$-9D_n^2-\left(K_{X_C}\cdot D_n\right)^2=9d^2-\left(\sum_{i=0}^n u_i\right)^2=\left(3d+\sum_{i=0}^n
u_i\right)\left(-K_{X_C}\cdot \tilde{C}^{\ck_C}\right),$$  where
the first one follows by Bézout's Theorem.

The equivalence between $(b)$ and $(c)$ follows from the
expression of the arithmetic genus of $\tilde{C}^{\ck_C}$,
$$p_a(\tilde{C}^{\ck_C})=1+\frac{1}{2}K_{X_C}\cdot
\tilde{C}^{\ck_C},$$ and the following fact: if $C$ were not
smooth in its affine part, then its geometric genus would be less
than $p_a(\tilde{C}^{\ck_C})$. \findemo\\

A direct consequence of this proposition is the following

\begin{cor}\label{www}
If $C$ is a curve of AMS type, then all the curves of the pencil
$\cp(C)$, except the non-reduced one, are also of AMS type.
\end{cor}

\begin{rem}\label{rectarat2}
{\rm If $C$ is a curve of AMS type, then $\tilde{H}^{\ck_C}$ is a
$(-1)$-curve of $X_C$. This fact is trivial from what is said in
Section \ref{ams}.}
\end{rem}

The following result provides a characterization of the proximity
graphs of the form $\gG({\mathcal K}_C)$, where $C$ is a curve of
AMS type.

\begin{pro}\label{sss}
Let $\gG$ be a proximity graph. Then, there exists a curve $C$ of
AMS type such that $\gG=\gG({\mathcal K}_C)$ if and only if there is
a sequence of integers $n_1,n_2,\ldots,n_r\geq 2$ such that
$\gG=\gG(n_1,n_2,\ldots,n_r)^+$.
\end{pro}

\noindent {\it Proof}. Assume the existence of a curve $C$ of AMS
type such that $\gG=\gG({\mathcal K}_C)$. The degree of $C$ is $n_1
n_2\cdots n_r$ by Corollary \ref{pernia}. Suppose that the
configuration ${\mathcal K}_C$ is $(p_0,p_1,\ldots,p_n)$ and let
${\mathcal C}$ be the configuration provided by the minimal embedded
resolution of $C$. Since the morphism $\pi_{{\mathcal K}_C}: X_C
\rightarrow \gp^2$ induces an embedded resolution of $C$, one has
that ${\mathcal C}=(p_0,p_1,\ldots,p_m)$ for some $m\leq n$. Then,
applying Bézout's Theorem to two curves of the pencil $\cp(C)$, the
following equality holds:
\begin{equation}\label{b}
(n_1 n_2\cdots n_r)^2=\sum_{i=0}^m u_i^2 +n-m,
\end{equation}
where $u_i$ denotes the multiplicity of the strict transform of $C$
at $p_i$. Taking into account that $\gG({\mathcal
C})=\gG(n_1,n_2,\ldots,n_r)^-$ and $u_i=\sum_{p_j\rightarrow p_i}
u_j$ for all $i=1,2,\ldots,m$, the multiplicities $u_i$ can be
easily computed in terms of $n_1,n_2,\ldots,n_r$ and, from the
equality $(\ref{b})$, it holds that $n-m=n_r$. Therefore,
$\gG({\mathcal K}_C)=\gG(n_1,n_2,\ldots,n_r)^+$.

Conversely, assume that $\gG=\gG(n_1,n_2,\ldots,n_r)^+$ for certain
integers $n_1,n_2,\ldots,n_r\geq 2$. By Proposition \ref{villa},
there exists a curve $C$ of AMS type such that the proximity graph
associated with its minimal embedded resolution is
$\gG(n_1,n_2,\ldots,n_r)^-$. Similar arguments to those used in the
above paragraph show that $\gG({\mathcal K}_C)=\gG$.\findemo

\begin{de}
{\rm For each curve $C$ having one place at infinity, the {\it
effective} (resp., {\it nef}) {\it semigroup} of $X_C$, denoted by
$NE(X_C)$ (resp., $P(X_C)$), is defined as the subsemigroup of
$\pic(X_C)$ generated by the classes of all effective (resp.,
numerically effective) divisors on $X_C$. }
\end{de}

Campillo, Piltant and Reguera described, in \cite{c-p-r1}, the
effective semigroup of $X_C$. They proved the following equality:
$$NE(X_C)=\gn [\tilde{H}^{\ck_C}]\oplus \bigoplus_{i=0}^n \gn
[\tilde{E}_i^{\ck_C}],$$ where $H$ denotes, as above, the line of
infinity. Moreover, as a consequence of \cite[Cor. 7 and Prop.
6]{c-p-r1} and Corollary \ref{www}, we have the following result:

\begin{pro}\label{sup}
If $C$ is a curve of AMS type, then $P(X_C)$ coincides with the
semigroup of classes in $\pic(X_C)$ of the form
$[\tilde{D}^{\ck_C}]$, where $D\hookrightarrow \gp^2$ is a
projective curve whose support does not contain the line of
infinity.
\end{pro}

Now, we will state and prove the announced $H^1$-vanishing result
for line bundles on $X_C$.

\begin{teo}\label{gordo}
Let $C$ be a curve of AMS type.
\begin{itemize}
\item[(a)] If $D$ is a numerically effective divisor on $X_C$,
then $h^1(X_C, \co_{X_C}(D))=0$.
\item[(b)] Let $D$ be an effective divisor on $X_C$ such that
$D\cdot \tilde{E}_i^{\ck_C}\geq 0$ for all $i=0,1,\ldots,n$ and
$D\cdot \tilde{E}_1^{\ck_C}\geq 1$ . Then, $h^1(X_C,
\co_{X_C}(D))=0$ if and only if $D\cdot \tilde{H}^{\ck_C}\geq -1$.

\end{itemize}

\end{teo}
\noindent {\it Proof}. Firstly, notice that the configuration
$\ck_C$ is P-sufficient, by Proposition \ref{rat}. In order to
prove $(a)$, we will reason by contradiction. So, we assume the
existence of a numerically effective divisor $D$ such that
$h^1(X_C, \co_{X_C}(D))>0$. Since $[D]$ is an effective class (by
Proposition \ref{sup}) we can apply Proposition \ref{2} and
\cite[Lem. II.7]{harb}, deducing that the complete linear system
$|D|$ has fixed part. Moreover, again by Proposition \ref{sup},
this fixed part has not exceptional components (that is, it has no
divisor $\tilde{E}_i^{\ck_C}$ as a component).

Now, the integral fixed components of $|D|$ have negative
self-intersection. Indeed, if we assume the existence of an
integral fixed component $R$ such that $R^2\geq 0$, then
$K_{X_C}\cdot R<0$ (by Proposition \ref{2}) and, applying the
Riemann-Roch Formula, we get
$$h^0(X_C,\co_{X_C}(R))\geq 1+(R^2-K_{X_C}\cdot R)/2\geq 2,$$
which is false, since $h^0(X_C,\co_{X_C}(R))=1$.

Finally Proposition \ref{sup} provides a contradiction, because
the unique non-exceptional integral curve on $X_C$ with negative
self-intersection is $\tilde{H}^{\ck_C}$.

To prove $(b)$, we consider a divisor $D$ satisfying the hypotheses.
First, we will assume the inequality $h^1(X_C, \co_{X_C}(D))>0$ and
we will show that this implies $D\cdot \tilde{H}^{\ck_C}\leq -2$. By
applying $(a)$ we have $D\cdot \tilde{H}^{\ck_C}\leq -1$ and,
therefore, it only remains to prove that the equality $D\cdot
\tilde{H}^{\ck_C}=-1$ leads to a contradiction. Using the hypotheses
and Remark \ref{rectarat2}, it can be deduced that
$D-\tilde{H}^{\ck_C}$ is a numerically effective divisor. A similar
reasoning to that given in the proof of $(a)$ shows that
$\tilde{H}^{\ck_C}$ is the unique possible integral fixed component
of the linear system $|D-\tilde{H}^{\ck_C}|$ and, then, it must be
fixed part free by Proposition \ref{sup}. So, we have a
decomposition $[D]=[\tilde{H}^{\ck_C}]+[T]$, where $T$ is an
effective divisor such that $|T|$ is fixed part free. Applying Part
$(a)$ to the divisor $T$, we deduce that $h^1(X_C,\co_{X_C}(T))=0$.
Taking into account this fact, Riemann-Roch Theorem and the equality
$h^0(X_C,\co_{X_C}(D))=h^0(X_C,\co_{X_C}(T))$, the following chain
of equalities and inequalities holds:
$$0<h^1(X_C,\co_{X_C}(D))=h^0(X_C,\co_{X_C}(D))-1-\frac{1}{2}(D^2-K_{X_C}\cdot
D)=$$
$$=h^0(X_C,\co_{X_C}(T))-1-\frac{1}{2}(D^2-K_{X_C}\cdot
D)=$$ $$=1+\frac{1}{2}(T^2-K_{X_C}\cdot
T)-1-\frac{1}{2}(D^2-K_{X_C}\cdot D)=$$
$$=\frac{1}{2}(-(\tilde{H}^{\ck_C})^2-2 \tilde{H}^{\ck_C}\cdot T+K_{X_C}\cdot \tilde{H}^{\ck_C}).$$
Hence, $(\tilde{H}^{\ck_C})^2-K_{X_C}\cdot
\tilde{H}^{\ck_C}<-2\tilde{H}^{\ck_C}\cdot T\leq 0$. But this is a
contradiction, since $\tilde{H}^{\ck_C}$ is a $(-1)$-curve.

It only remains to prove that, if $D\cdot \tilde{H}^{\ck_C}\leq -2$,
then $h^1(X_C, \co_{X_C}(D))>0$. But the inequality $D\cdot
\tilde{H}^{\ck_C}\leq -2$ implies that the $(-1)$-curve
$\tilde{H}^{\ck_C}$ is a multiple fixed component of the linear
system $|D|$, and it is easy to see that this fact implies that
$h^1(X_C, \co_{X_C}(D))$ is positive (by a similar reasoning to that
given in \cite[pag. 197]{mir}).\findemo

\begin{rem}\label{dim}
{\rm  With the hypotheses of Theorem \ref{gordo}, Part $(a)$
allows us to determine the dimension of all complete linear
systems on $X_C$. Indeed, let $D$ be a divisor on $X_C$ and
consider the set $S=\{[\tilde{H}^{\ck_C}],
[\tilde{E}_0^{\ck_C}],[\tilde{E}_1^{\ck_C}],\ldots,[\tilde{E}_n^{\ck_C}]\}\subseteq
\pic(X_C)$. If for some $F\in S$ we have $D\cdot F<0$, then it is
obvious that $h^0(X_C,\co_{X_C}(D))=h^0(X_C,\co_{X_C}(D-F))$.
Therefore, we can perform the following process: check $D\cdot F$
for each $F\in S$, replace $D$ by $D-F$ whenever $D\cdot F<0$ and
continue with the new $D$. The process ends when it gives rise to
a divisor $D'$ such that either it is obviously not effective
(because either $D'\cdot L^{\ck_C}<0$ or
$D'\cdot(L^{\ck_C}-E_1^{\ck_C})<0$) or $D'$ is numerically
effective. Since $h^0(X_C,\co_{X_C}(D))=h^0(X_C,\co_{X_C}(D'))$,
in the first case the linear system $|D|$ is empty and in the
second case its dimension is  $(D'^2-K_{X_C}\cdot D')/2$, by Part
$(a)$ of Theorem \ref{gordo}.

}
\end{rem}

\section{Linear systems of curves through generic points on
$\gp^2$}\label{s4}

In this section we will use Theorem \ref{gordo} and a
specialization process to deduce some results about the dimension
of linear systems of curves passing through a finite set of points
of the plane in generic position.

\subsection{Special linear systems and the Harbourne-Hirschowitz
Conjecture}\label{s41}

Given a projective curve $C$ of $\gp^2$, we will say that $C$ {\it
goes through} a weighted configuration
$(\ck=(p_i)_{i=0}^n,\bold{m}=(m_i)_{i=0}^n)$ if and only if the
divisor $C^{\ck}-\sum_{i=0}^n m_i E_i^{\ck}$ is effective. For any
degree $d$, denote by $\cL_d(\ck,\bold{m})$ the set of projective
curves on $\gp^2$ of degree $d$ going through $(\ck,\bold{m})$.
This is a linear system of $\gp^2$ that is projectively isomorphic
to the complete linear system $|D_{d,\ck,\bold{m}}|$ of $Z_{\ck}$,
$D_{d,\ck,\bold{m}}$ being the divisor $d L^{\ck}-\sum_{i=0}^n m_i
E_i^{\ck}$. From Riemann-Roch Theorem one gets
$$\dim \cL_d(\ck,\bold{m})-h^1(\cL_d(\ck,
\bold{m}))=\frac{d(d+3)}{2}-\sum_{i=0}^n \frac{m_i(m_i+1)}{2},$$
where $\dim \cL_d(\ck,\bold{m})$ is the dimension of
$\cL_d(\ck,\bold{m})$ as projective space and $h^1(\cL_d(\ck,
\bold{m})):=h^1(Z_{\ck},\co_{Z_{\ck}}(D_{d,\ck,\bold{m}}))$ will
be called the {\it superabundance} of $\cL_d(\ck,\bold{m})$. The
independence of the linear conditions imposed by the weighted
configuration $(\ck,\bold{m})$ is equivalent to the vanishing of
this superabundance.

If $\ck$ is a configuration whose points lie all in $\gp^2$, the
dimension and the superabundance of a linear system
$\cL_d(\ck,\bold{m})$ depend on the position of the points of $\ck$,
and they reach their minimal values for a generic set of points. We
will denote by $\ck_0(n)$ a configuration consisting of $n+1$ closed
points of $\gp^2$ in generic position. For each integer $d\geq 1$
and for each system of multiplicities
$\bold{m}=(m_0,m_1,\ldots,m_n)$ we will denote by $\cL_d(\bold{m})$
the linear system $\cL_d(\ck_0(n),\bold{m})$. Also, we define the
{\it expected dimension} of $\cL_d(\bold{m})$ to be the following
number:
$$\edim
\cL_d(\bold{m}):=\max \left\{\frac{d(d+3)}{2}-\sum_{i=0}^n
\frac{m_i(m_i+1)}{2},-1\right\}.$$

\begin{de}
{\rm We will say that a linear system $\cL_d(\bold{m})$ is {\it
special} if and only if $\dim \cL_d(\bold{m})>\edim
\cL_d(\bold{m})$, that is, $\cL_d(\bold{m})$ is non-empty and the
superabundance $h^1(\cL_d(\bold{m}))$ is positive.}
\end{de}

Given a positive integer $d$ and a system of multiplicities
$\bold{m}=(m_i)_{i=0}^n$, it is easy to prove that, if there exists
a curve $C$ on $\gp^2$ such that its strict transform on
$Z_{\ck_0(n)}$ is a $(-1)$-curve and $D_{d,\ck_0(n),\bold{m}}\cdot
\tilde{C}^{\ck_0(n)}\leq -2$, then the linear system ${\mathcal
L}_d(\bold{m})$ is special (see, for instance,  \cite[pag.
197]{mir}). One of the equivalent statements of the
Harbourne-Hirschowitz Conjecture
is just the converse assertion:\\

\noindent {\bf Conjecture}. \;(Harbourne-Hirschowitz) If a linear
system $\cL_d(\bold{m})$ is special, then there exists a curve $C$
on $\gp^2$ such that its strict transform on $Z_{\ck_0(n)}$ is a
$(-1)$-curve and $D_{d,\ck_0(n),\bold{m}}\cdot
\tilde{C}^{\ck_0(n)}\leq -2$.

\subsection{A non-speciality result and some consequences}\label{s42}

For each positive integer $n$, there exists a variety $Y_{n}$
whose points are naturally identified with the configurations over
$\gp^2$ with $n+1$ points. These varieties, known as {\it iterated
blowing-ups}, were introduced by Kleiman in \cite{kle1} and
\cite{kle2} and they have also been treated in \cite{harb2},
\cite{roe} and \cite{fdb2} (see also \cite{roe2}). There is a
family of projective morphisms $Y_{n+1}\rightarrow Y_{n}$ and
relative divisors $F_{-1},F_0,F_1,\ldots,F_n$ on $Y_{n+1}$ such
that the fiber over a given configuration
$\ck=(p_0,p_1,\ldots,p_n)$ (viewed as a point of $Y_{n}$) is
isomorphic to the surface $Z_{\ck}$ obtained by blowing-up the
points in $\ck$ and, if $i\geq 0$ (resp., $i=-1$), the restriction
of $F_i$ to this fiber corresponds to the total transform
$E_i^{\ck}$ of the exceptional divisor appearing in the blowing-up
centered at $p_i$ (resp., the total transform of a general line of
$\gp^2$).

For each positive integer $d$ and for each sequence of
multiplicities ${\bold m}=(m_0,m_1,\ldots,m_n)$ we apply the
Semicontinuity Theorem \cite[III, 12.8]{har} to the invertible
sheaf $\co_{Y_{n+1}}(dF_{-1}-m_0F_0-m_1F_1-\ldots-m_nF_n)$,
obtaining that the functions $Y_n\rightarrow \gz$ given by
\begin{equation}\label{semic}
{\ck} \mapsto h^i(Z_{\ck}, \co_{Z_{\ck}}(D_{d,\ck, \bold{m}})),
\end{equation}
for $i\in \{0,1\}$, are upper-semicontinuous.

For each proximity graph $\gG$ with $n+1$ vertices, we define
$U(\gG)$ as the subset of $Y_n$ containing exactly the
configurations $\ck$ whose proximity graph is $\gG$. This is an
irreducible smooth locally closed subvariety (\cite{roe} and
\cite{fdb2}). As a consequence of the upper-semicontinuity of the
functions given in (\ref{semic}), for any positive integer $d$ and
any system of multiplicities $\bold{m}=(m_0,m_1,\ldots,m_n)$, the
dimension and the superabundance of the linear systems ${\mathcal
L}_d(\ck,\bold{m})$, for $\ck$ varying in $U(\gG)$, achieve the
minimum value in a dense open subset of $U(\gG)$.

\begin{de}
{\rm We will say that a weighted configuration $(\ck,
\bold{m}=(m_i)_{i=0}^n)$ (resp., a weighted proximity graph
$(\gG,\bold{m})$) is {\it consistent} if all the excesses
$\rho_j(\ck,\bold{m})$ (resp., $\rho_j(\gG,\bold{m})$)) are
non-negative (see Section \ref{s1} for the definition of excesses).
In this case, and provided that $n\geq 1$, we associate with $(\ck,
\bold{m})$ an integer, denoted by $\epsilon(\ck, \bold{m})$ (or
$\epsilon(\gG(\ck), \bold{m})$, since it depends only on the
weighted proximity graph) and defined to be either $1$, if
$\rho_1(\ck, \bold{m})\geq 1$, or $0$, if $\rho_1(\ck, \bold{m})=0$.
}
\end{de}

Given a weighted proximity graph
$(\gG=\gG(\ck),\bold{m}=(m_i)_{i=0}^n)$,  it is possible to obtain a
unique system of multiplicities, which will be denoted by
$\bold{m}^{\gG}=(m_0^{\gG},m_1^{\gG},\ldots,m_n^{\gG})$, such that
$(\gG, \bold{m}^{\gG})$ is consistent and the ideal sheaves given by
${\pi_{\ck}}_*\co_{Z_{\ck}}(-\sum_{i=0}^n m_i E_i^{\ck})$ and
${\pi_{\ck}}_*\co_{Z_{\ck}}(-\sum_{i=0}^n m_i^{\gG} E_i^{\ck})$
coincide. So, there exists a canonical bijection between the linear
systems $|D_{d,\ck,\bold{m}}|$ and $|D_{d,\ck,\bold{m}^{\gG}}|$ for
all integers $d\geq 1$. The procedure used to obtain
$\bold{m}^{\gG}$ is called {\it unloading} \cite[4.6]{cas} and it
depends only on the proximity graph $\gG$, and not on a special
election of the configuration $\ck$ associated with $\gG$. In each
step of the unloading procedure ({\it unloading step}) one must
detect a point $p_i$ of $\ck$ such that its associated excess
$\rho_i(\ck,\bold{m})$ is negative; then, one replaces the system of
multiplicities $\bold{m}$ by the system
$\bold{m}'=(m_0',m_1',\ldots,m_n')$ where $m_i'=m_i+1$, $m_j'=m_j-1$
for those indexes $j$ such that $p_j$ is proximate to $p_i$, and
$m_j'=m_j$ otherwise (if some multiplicity in $\bold{m}'$ is
negative, it must be replaced by $0$). Now, we must perform another
unloading step from the new system $\bold{m}'$, and so on. A finite
number of unloading steps lead to the desired system of
multiplicities $\bold{m}^{\gG}$. An unloading step applied to a
point $p_i$ whose associated excess equals $-1$ is called {\it
tame}. Tame unloadings will be very useful for us, since they
preserve independence of conditions, that is, if $\bold{m}'$ is
obtained from $\bold{m}$ performing a tame unloading step, then
$h^1(Z_{\ck},\co_{Z_{\ck}}(D_{d,\ck,\bold{m}}))=h^1(Z_{\ck},\co_{Z_{\ck}}(D_{d,\ck,\bold{m}'}))$
for all positive integer $d$ (this fact can be easily deduced from
\cite[4.7.1]{cas} and \cite[4.7.3]{cas}).

\begin{de}
{\rm  We will say that a weighted configuration $(\ck, \bold{m})$
(resp., a weighted proximity graph $(\gG, \bold{m})$) is {\it
almost consistent} if either it is consistent or there exists a
sequence of tame unloading steps leading from $\bold{m}$ to
$\bold{m}^{\gG(\ck)}$ (resp., $\bold{m}^{\gG}$). }
\end{de}

The following theorem gives a sufficient condition for the
non-speciality of a linear system ${\mathcal L}_d(\bold{m})$, when
$\bold{m}$ is a system of multiplicities such that
$(\ck_C,\bold{m})$ is almost consistent, $C$ being a curve of AMS
type. Moreover, when this weighted configuration is consistent and
the excess $\rho_1(\ck_C,\bold{m})$ is positive, it provides a
characterization of such non-special linear systems which are not
empty.

\begin{teo}\label{super}
Let $d$ be a positive integer and $\bold{m}=(m_0,m_1,\ldots,m_n)$ a
system of multiplicities, with $n\geq 1$. Assume the existence of a
curve $C$ of AMS type such that $(\ck_C,\bold{m})$ is almost
consistent. Then, the linear system ${\mathcal L}_d(\bold{m})$ is
non-special whenever $d\geq m_0^{\gG}+m_1^{\gG}-\epsilon(\gG,
\bold{m}^{\gG})$, where $\gG:=\gG(\ck_C)$. Moreover, if ${\mathcal
L}_d(\bold{m})$ is not empty, $(\ck_C,\bold{m})$ is consistent and
$\rho_1(\gG, \bold{m})\geq 1$, then the following equivalence holds:
${\mathcal L}_d(\bold{m})$ is non-special if and only if $d\geq
m_0+m_1-1$.
\end{teo}
\noindent {\it Proof}. Set $\ck_C=(p_0,p_1,\ldots,p_n)$ (adding null
multiplicities to $\bold{m}$, if it is necessary, we can assume that
the cardinality of $\ck_C$ is $n+1$). In order to prove the first
assertion of the statement, we will reason by contradiction. So,
assume that $d\geq m_0^{\gG}+m_1^{\gG}-\epsilon(\gG,
\bold{m}^{\gG})$ and ${\mathcal L}_d(\bold{m})$ is special. The
subset $U(\gG(\ck_0(n)))$ is dense in $Y_n$ (see \cite{kle1}) and
then, as a consequence of the upper-semicontinuity of the functions
given in (\ref{semic}), the following inequalities hold: $\dim
{\mathcal L}_d(\bold{m})\leq
h^0(X_C,\co_{X_C}(D_{d,\ck_C,\bold{m}}))-1$ and $h^1({\mathcal
L}_d(\bold{m}))\leq h^1(X_C,\co_{X_C}(D_{d,\ck_C,\bold{m}}))$. Thus,
the complete linear system on $X_C$ given by
$|D_{d,\ck_C,\bold{m}}|$ is not empty and
$h^1(X_C,\co_{X_C}(D_{d,\ck_C,\bold{m}}))$ is positive.  But, since
$(\ck_C,\bold{m})$ is almost consistent, we have
$$h^1(X_C,\co_{X_C}(D_{d,\ck_C,\bold{m}^{\gG}}))=h^1(X_C,\co_{X_C}(D_{d,\ck_C,\bold{m}})).$$
The consistency of $(\ck_C,\bold{m}^{\gG})$ implies that
$D_{d,\ck_C,\bold{m}^{\gG}}\cdot \tilde{E}_i^{\ck_C}\geq 0$ for
all $i=0,1,\ldots,n$ and, therefore, we can apply Theorem
\ref{gordo} to deduce the inequality
$D_{d,\ck_C,\bold{m}^{\gG}}\cdot \tilde{H}^{\ck_C}\leq
-1-\epsilon(\gG, \bold{m}^{\gG})$. But, taking into account that
$\tilde{H}^{\ck_C}$ is a $(-1)$-curve, this is equivalent to the
condition $d\leq m_0^{\gG}+m_1^{\gG}-1-\epsilon(\gG,
\bold{m}^{\gG})$, a contradiction.

For the last assertion, it only remains to prove that, assuming the
consistency of $(\gG,\bold{m})$ and the inequality $\rho_1(\gG,
\bold{m})\geq 1$, the non-speciality of the linear system ${\mathcal
L}_d(\bold{m})$ implies the inequality $d\geq m_1+m_2-1$. We will
reason by contradiction. So, assume that ${\mathcal L}_d(\bold{m})$
is non-special and $d\leq m_1+m_2-2$. Again taking into account that
$\tilde{H}^{\ck_C}$ is a $(-1)$-curve, we have
$D_{d,\ck_C,\bold{m}}\cdot \tilde{H}^{\ck_C}\leq -2$. If $N$ denotes
the line of $\gp^2$ joining the two first points of the
configuration $\ck_0(n)$, then
$$D_{d,\ck_0(n),\bold{m}}\cdot \tilde{N}^{\ck_0(n)}=d-m_0-m_1=D_{d,\ck_C,\bold{m}}\cdot
\tilde{H}^{\ck_C}\leq -2,$$ which is a contradiction with the
non-speciality of  ${\mathcal L}_d(\bold{m})$.\findemo\\

As a consequence of Proposition \ref{sss}, there is a bijection
between the set of ordered sequences $(n_1,n_2,\ldots,n_r)\in
(\gn\setminus \{0,1\})^r$ (with $r\in \gn\setminus \{0\}$) and the
set of proximity graphs of the form $\gG(\ck_C)$, $C$ being a
curve of AMS type. Taking this fact into account, we obtain the
following reformulation of Theorem \ref{super}, expressed in
purely arithmetical terms:

\begin{cor}\label{super2}
Let $d$ be a positive integer and $\bold{m}=(m_0,m_1,\ldots,m_n)$ a
system of multiplicities with $n\geq 1$. Assume that there exist
integers $n_1,n_2,\ldots,n_r\geq 2$  such that the weighted
proximity graph $(\gG=\gG(n_1,n_2,\ldots,n_r)^+,\bold{m})$ is almost
consistent. Then, the linear system ${\mathcal L}_d(\bold{m})$ is
non-special whenever $d\geq
m_0^{\gG}+m_1^{\gG}-\epsilon(\gG,\bold{m})$. Moreover, if ${\mathcal
L}_d(\bold{m})$ is not empty, $(\gG,\bold{m})$ is consistent and
$\rho_1(\gG,\bold{m})\geq 1$, then the following equivalence holds:
${\mathcal L}_d(\bold{m})$ is non-special if and only if $d\geq
m_0+m_1-1$.
\end{cor}

Next, we will give two examples by applying Corollary \ref{super2}
to two specific sequences of integers.

\begin{exa}\label{ex1}
{\rm Let $n\geq 2$ be an integer and consider the proximity graph
$\gG=\gG(t+1)^+$, where $t:=\lfloor n/2 \rfloor$. The number of
points $s+1$ of whichever configuration whose proximity graph be
$\gG$ is $n+1$ (resp., $n+2$) if $n$ is odd (resp., if $n$ is even)
and, moreover, the complete list of proximity relations among the
points of such a configuration $(p_0,p_1,\ldots,p_s)$ are the
following: $p_i\rightarrow p_{i-1}$ for all $i=1,2,\ldots,s$ and
$p_j\rightarrow p_0$ for all $j=2,3,\ldots,t$. Then, applying
Corollary \ref{super2} to the graph $\gG$, we get the following
result:

Let $\bold{m}=(m_0,m_1,\ldots,m_n)$ be a system of multiplicities
such that $m_1\geq m_2\geq \ldots \geq m_n$ and $m_0\geq
\sum_{i=1}^t m_i$. A linear system ${\mathcal L}_d(\bold{m})$ is
non-special whenever $d\geq m_0+m_1-\epsilon$, where
$\epsilon=\min\{1,m_1-m_2\}$. Moreover, if ${\mathcal
L}_d(\bold{m})$ is not empty and $m_1>m_2$, then ${\mathcal
L}_d(\bold{m})$ is non-special if and only if $d\geq m_0+m_1-1$. }
\end{exa}

\begin{exa}\label{ex2}
{\rm Let $k\geq 2$ be an integer. If $(p_0,p_1,\ldots,p_s)$ is
whichever configuration whose proximity graph is
$\gG(2,2,\ldots,2)^+$ (where the number $2$ appears $k$ times), one
gets that $s=3k$ and the proximity relations among the points are
the following: $p_i\rightarrow p_{i-1}$ for all $i=1,2,\ldots,s$ and
$p_{3j+1}\rightarrow p_{3j-1}$ for all $j=1,2,\ldots,k-1$. Applying
again Corollary \ref{super2} to this graph one gets the following
result:

Let $\bold{m}=(m_0,m_1,\ldots,m_{3k})$ be a system of multiplicities
such that $m_0\geq m_1\geq \ldots \geq m_{3k}$ and $m_{3i-1}\geq
m_{3i}+m_{3i+1}$ for all $i=1,2,\ldots,k-1$. Then, a linear system
${\mathcal L}_d(\bold{m})$ is non-special if $d\geq m_0+m_1$. If, in
addition, ${\mathcal L}_d(\bold{m})$ is not empty and $m_1>m_2$,
then ${\mathcal L}_d(\bold{m})$ is non-special if and only if $d\geq
m_0+m_1-1$. }
\end{exa}

The following direct consequence of Theorem \ref{super} exhibits a
wide range of cases in which the Harbourne-Hirschowitz Conjecture
is satisfied.

\begin{cor}\label{conj}
Let $\bold{m}=(m_0,m_1,\ldots,m_n)$ be a system of multiplicities
such that $(\ck_C, \bold{m})$ is consistent and
$\rho_1(\ck_C,\bold{m})\geq 1$, $C$ being a curve of AMS type.
Denote $\ck_0(n)=(p_0,p_1,\ldots,p_n)$ and set $N$ the line joining
$p_0$ and $p_1$.  If a linear system of the form ${\mathcal
L}_d(\bold{m})$ is special, then $D_{d,\ck_0(n),\bold{m}}\cdot
\tilde{N}^{\ck_0(n)}\leq -2$.
\end{cor}

\begin{rem}\label{nequivalence}
{\rm For a fixed positive integer $n$, let us denote by ${\mathcal
S}_n$ the set of proximity graphs of the form $\gG({\mathcal K}_C)$,
where $C$ is a curve of AMS type, whose number of vertices is
greater than or equal to $n+1$. Each proximity graph
$\gG=\gG({\mathcal K}_C)\in {\mathcal S}_n$ provides, by Corollary
\ref{conj}, an infinite family of systems of multiplicities
$\bold{m}=(m_i)_{i=0}^n$ for which the Harbourne-Hirschowitz
Conjecture is true. This family is given by the non-negative integer
solutions of the following system of linear inequalities in
$m_0,m_1,\ldots,m_n$:
$$\left\{ \begin{array}{l}
 m_i  - \sum\limits_{j\leq n;\; p_j  \to p_i } {m_j  \ge 0,\quad i = 0,2,3,4, \ldots ,n}  \\
 m_1  - \sum\limits_{j\leq n;\; p_j  \to p_1 } {m_j  \ge 1}  \\
 \end{array} \right.$$
where $\ck_C=(p_0,p_1,\ldots,p_n)$. Two graphs $\gG$ and $\gG'$ in
${\mathcal S}_n$ give rise to the same system of inequalities if
the proximity relations involving the first $n+1$ vertices are the
same for both graphs; in this case, we will say that $\gG$ and
$\gG'$ are {\it $n$-equivalent}. Taking into account Proposition
\ref{sss}, a complete system of representants of the quotient set
of ${\mathcal S}_n$ by this equivalence relation is given by the
proximity graphs $\gG(n_1,n_2,\ldots,n_r)^+$ (with $n_i\geq 2$ for
all $i$) such that either $r=1$ and $(n+1)/2\leq n_1\leq n+1$, or
$r>1$, $t:=n-2\sum_{i=1}^{r-1} n_i+r>0$ and $t/2\leq n_r\leq t$.
Thus, for a fixed positive integer $n$, the set of distinct
systems of linear inequalities in $n+1$ variables provided by
Corollary \ref{conj} (each of them satisfying that the
Harbourne-Hirschowitz Conjecture is true for the solutions) is
finite. In fact, they are in one-to-one correspondence with the
$n$-equivalence classes. }
\end{rem}

\begin{de}
{\rm We define the {\it regularity} of a system of multiplicities
$\bold{m}=(m_0,m_1,\ldots,m_n)$ as the minimum integer $d$ such
that $(\ck_0(n),\bold{m})$ imposes independent conditions to the
curves of degree $d$, and we will denote it by $\tau(\bold{m})$. }
\end{de}

The following result is another consequence of Theorem \ref{super}
and it allows to compute the exact value of the regularity of a wide
range of systems of multiplicities:

\begin{cor}\label{regu}
Let $\bold{m}=(m_0,m_1,\ldots,m_n)$ be a system of multiplicities,
with $n\geq 1$. Assume the existence of a curve $C$ of AMS type such
that $(\ck_C,\bold{m})$ is consistent and
$\rho_1(\ck_C,\bold{m})\geq 1$. Then, $\tau(\bold{m})=m_0+m_1-1$ if
$(m_0+m_1-1)(m_0+m_1+2)-\sum_{i=0}^n m_i(m_i+1)\geq -2$, and
$\tau(\bold{m})=m_0+m_1$ otherwise.
\end{cor}
\noindent {\it Proof}. First, we will show that
$[D_{d,\ck_0(n),\bold{m}}]$ is an effective class, where
$d=m_0+m_1$. Indeed, the class $[D_{d,\ck_C,\bold{m}}]$ is
numerically effective, due to the hypotheses and the fact that
$\tilde{H}^{{\mathcal K}_C}$ is a $(-1)$-curve. So, it is an
effective class of $\pic(X_C)$ by Proposition \ref{sup}, and
$h^1({\mathcal
L}_d(\ck_C,\bold{m}))=h^1(X_C,\co_{X_C}(D_{d,\ck_C,\bold{m}}))=0$,
by Theorem \ref{gordo}. Therefore, we get
$$\edim {\mathcal L}_d(\bold{m})=\dim {\mathcal L}_d(\ck_C,\bold{m})\geq
0.$$ Hence, ${\mathcal L}_d(\bold{m})$ is non-empty and then, by
Theorem \ref{super}, $\tau(\bold{m})\leq d$.

If $(d-1)(d+2)-\sum_{i=0}^n m_i(m_i+1)\geq -2$, then the inequality
$\edim {\mathcal L}_{d-1}(\bold{m})\geq -1$ holds. So, the
superabundance $h^1({\mathcal L}_{d-1}(\bold{m}))$ must be zero, in
virtue of Theorem \ref{super}, and then $\tau(\bold{m})=d-1$ by
\cite[2]{hir}.

Finally, if $(d-1)(d+2)-\sum_{i=0}^n m_i(m_i+1)< -2$, then $$\dim
{\mathcal L}_{d-1}(\bold{m})-h^1({\mathcal L}_{d-1}(\bold{m}))<-1$$
and this implies that the superabundace  $h^1({\mathcal
L}_{d-1}(\bold{m}))$ is positive. Therefore, in this case,
$\tau(\bold{m})=d$.\findemo

\subsection{Bounding the regularity}\label{s43}

In \cite{roe3} it is described an algorithm, based on the
unloading method, which provides an upper bound of the regularity
of whichever system of multiplicities
$\bold{m}=(m_0,m_1,\ldots,m_n)$. Although only the case of
homogeneous multiplicities is explicitly treated (i.e.,
$m_0=m_1=\cdots=m_n$) this algorithm can be adapted without
difficulty to the case of arbitrary multiplicities. In this
section we introduce a generalization of this algorithm, based on
our results in Section \ref{infinity}.

Let $\bold{m}=(m_0,m_1,\ldots,m_n)$ be a sequence of
multiplicities (with $n\geq 1$) such that $m_0\geq m_1\geq \cdots
\geq m_{n}$. Take a sequence of integers $(n_1, n_2, \ldots,n_r)$
such that $n_i\geq 2$ for all $i=1,2,\ldots,r$  and $n+1$ is not
greater than the number of vertices of the graph $\gG:=\gG(n_1,
n_2,\ldots,n_r)^+$. By completing with zero multiplicities, if it
is necessary, we will assume that this number of vertices
coincides with $n+1$. Denote by $\gG_i$ the proximity graph
obtained from $\gG$ by deleting all curved-dotted edges involving
some vertex with label greater than $i$. Let $(i_1,
i_2,\ldots,i_w)$ be an increasing sequence of integers such that
$\gG_{i_1}, \gG_{i_2},\ldots, \gG_{i_w}$ are the distinct elements
of the set $\{\gG_i \mid 1\leq i < n\}$ and let $a$ be the maximum
integer such that $0\leq a< n$ and, if $(p_0,p_1,\ldots,p_n)$
denotes a configuration with proximity graph $\gG$, the
cardinality of the set $\{ j\mid a\leq j\leq n, \;\;p_j\rightarrow
p_a\;\;\mbox{and}\;\; m_j>0\}$ is greater than $1$ (if that
integer does not exist, we will take $a=0$).

\begin{exa}
{\rm Let $\gG$ be the proximity graph $\gG(s)^+$, where $s>1$ is
an integer. The above described sequence of graphs $\gG_{i_1},
\gG_{i_2},\ldots, \gG_{i_w}$ is, in this case, the sequence
$\gG^1, \gG^2,\ldots,\gG^{s-1}$ where $\gG^1$ denotes the
proximity graph of a chain of $2s$ free points and, for each
$k=2,3,\ldots,s-1$, $\gG^k$ stands for the graph depicted in
Figure 3 of page \pageref{figure2} taking $n=2s-1$; the integer
$a$ is $0$.}
\end{exa}

Set $\bold{m}_1:=\bold{m}$, which is consistent for $\gG_{i_1}$,
and define recursively the systems of multiplicities
$\bold{m}_2,\bold{m}_3,\ldots,\bold{m}_w$ as follows. Suppose we
have defined $\bold{m}_k$ and perform the following two-steps
algorithm applied to $\bold{v}:=\bold{m}_k$, which will give rise
to $\bold{m}_{k+1}$:\\

\noindent {\it Step 1}. If $(\gG_{i_{k+1}},\bold{v})$ is
consistent, then define $\bold{m}_{k+1}:=\bold{v}$. Otherwise,
there exists a unique $j$ such that the excess
$\rho_j(\gG_{i_{k+1}},\bold{v})$ is negative. In this case, if
$j=a$ and this excess equals $-1$, define also
$\bold{m}_{k+1}:=\bold{v}$; else, perform an unloading step to
$(\gG_{i_k},\bold{v})$ at the vertex which corresponds to that
excess, replace $\bold{v}$ by the obtained new system of
multiplicities and go to
Step 2.\\

\noindent {\it Step 2}. Replace $\bold{v}$ by
$\bold{v}^{\gG_{i_k}}$ and return to Step 1.\\

Once we have computed $\bold{m}_w$,  we must consider the system
of multiplicities $\bold{m}':=\bold{m}_w^{\gG}=(m_0',m_1',\ldots,
m_n')$. Let $C$ be a curve having one place at infinity whose
associated proximity graph is $\gG_{i_w}=\gG$. Notice that
$h^1(X_C,\co_{X_C}(D_{m_0'+m_1', \ck_C, \bold{m}'}))=0$ (by
Theorem \ref{gordo}), since $D_{m_0'+m_1', \ck_C, \bold{m}'}$ is a
numerically effective divisor of $X_C$. We will compute the
successive dimensions $h^1(X_C,\co_{X_C}(D_{m_0'+m_1'-j, \ck_C,
\bold{m}'}))$ for $j=1,2,\ldots$ (using the process described in
Remark \ref{dim}) until finding the minimum $j$ such that the
mentioned dimension is positive. Finally, we will define
$\beta(\bold{m}):=m_0'+m_1'-j+1$. Note that this process is
independent of the chosen curve $C$; in fact, it only depends on
the proximity graph $\gG$.

Now, we will justify that the obtained value $\beta({\bold m})$ is
an upper bound of the regularity $\tau(\bold{m})$. We start with a
lemma whose proof is an adaptation of that of \cite[Lem.
2.1]{roe3} and we will omit it.

\begin{lem}\label{a}
Let $d$ be a positive integer, $\ck=\{p_0,p_1,\ldots, p_n\}$ a
configuration and $\bold{m}=(m_i)_{i=0}^n$ a system of
multiplicities. Let $i\in \{0,1,\ldots,n\}$ be such that
$\rho_i(\ck,\bold{m})\geq -1$ and let
$\bold{m}'=(m_0',m_1',\ldots,m_n')$ be the sequence of
multiplicities obtained from $\bold{m}$ by performing an unloading
step at the point $p_i$. Then, $h^1({\mathcal L}_d(\ck,\bold{m}))=0$
whenever $h^1 ({\mathcal L}_d(\ck,\bold{m}'))=0$.
\end{lem}

\begin{teo}\label{just}
Let $\bold{m}=(m_0,m_1,\ldots,m_n)$ be a sequence of
multiplicities, let $(n_1,n_2,\ldots,n_r)$ be a sequence of
integers such that $n_i\geq 2$ for all $i=1,2,\ldots,r$ and set
$\beta(\bold{m})$ defined as above. Then, $\beta(\bold{m})$ is an
upper bound of $\tau(\bold{m})$.
\end{teo}
\noindent {\it Proof}. Recall the notations of Section \ref{s42}.
Taking into account that the proximity graphs
$\gG_{i_1},\gG_{i_2},\ldots,\gG_{i_w}$ are associated with chain
configurations and the matrix ${\bf P}_{\gG_{i_{k-1}}}^{-1}\cdot
{\bf P}_{\gG_{i_k}}$ has no negative entries for each
$k=2,3,\ldots,w$, it can be deduced, from \cite{roe}, the
existence of a chain of inclusions
\begin{equation}\label{seq}
U(\gG)=U(\gG_{i_w})\subseteq \overline{U(\gG_{i_{w-1}})}\subseteq
\ldots \subseteq \overline{U(\gG_{i_1})}.
\end{equation}

Let $d= \beta(\bold{m})$ and, for each system of multiplicities
$\bold{v}$, set $h^1(d,\gG_{i_k},\bold{v})$ the minimum of the
superabundances $h^1({\mathcal L}_d(\ck,\bold{v}))$ when $\ck$
varies in $U(\gG_{i_k})$.

Take a plane curve $C$ of AMS type such that $\gG({\mathcal
K}_C)=\gG(n_1,n_2,\ldots,n_r)^+$. From the above description of the
algorithm, it follows that $$h^1({\mathcal
L}_d(\ck_C,\bold{m}'))=h^1(X_C,\co_{X_C}(D_{d,\ck_C,\bold{m}'}))=0$$
and, since the weighted proximity graph $(\gG,\bold{m}')$ is
obtained from $(\gG,\bold{m}_w)$ by tame unloading steps, we get
that the integer $h^1({\mathcal L}_d(\ck_C,\bold{m}_w))$ vanishes
and, hence, $h^1(d,\gG_{i_w},\bold {m}_w)=0$.

Finally, for $2\leq k \leq w$, we will show that the vanishing of
$h^1(d,\gG_{i_k},\bold{m}_k)$ implies that of
$h^1(d,\gG_{i_{k-1}},\bold{m}_{k-1})$. In order to prove this
assertion observe firstly that, if we assume that
$h^1(d,\gG_{i_k},\bold{m}_k)=0$, then
$h^1(d,\gG_{i_{k-1}},\bold{m}_k)=0$ by (\ref{seq}) and the
upper-semicontinuity of the functions given in (\ref{semic}). Choose
a configuration $\ck\in U(\gG_{i_{k-1}})$ such that $h^1({\mathcal
L}_d(\ck, \bold{m}_k))=0$.  It is not hard to see that the unloading
procedure of the Step 2 of the algorithm to obtain the sequence
$\bold{m}_1,\bold{m}_2,\ldots,\bold{m}_w$ can be performed by means
of tame unloading steps. From this fact and Lemma \ref{a}, the
equality $h^1({\mathcal L}_d(\ck, \bold{m}_{k-1}))=0$ is obtained
and, therefore, $h^1(d,\gG_{i_{k-1}},\bold{m}_{k-1})=0$.

Now,  it follows, by induction, that
$h^1(d,\gG_{i_1},\bold{m})=0$. Finally, using again semicontinuity
and taking into account the density of $U(\gG(\ck_0(n)))$ in
$Y_n$, we get $h^1(\cL_d(\bold{m}))=0$. Hence, $d$ is an upper
bound of $\tau(\bold{m})$.\findemo\\

We conclude the paper with some remarks on the above described
algorithmic bound.

First observe that, given a system of multiplicities $\bold{m}$,
there is a bound $\beta(\bold{m})$ for each election of a
proximity graph $\gG(n_1,n_2,\ldots,n_r)^+$ with, at least, $n+1$
vertices (that is, $2\sum_{i=1}^r n_i-r\geq n$). It is clear that
$n$-equivalent proximity graphs give rise to the same bound (see
Remark \ref{nequivalence}). Thus, one can apply the algorithm to
all the graphs $\gG(n_1,n_2,\ldots,n_r)^+$ such that either $r=1$
and $(n+1)/2\leq n_1\leq n+1$, or $r>1$, $t:=n-2\sum_{i=1}^{r-1}
n_i+r>0$ and $t/2\leq n_r\leq t$, and then pick the best bound.

We will show that the algorithm given by Roé in \cite{roe3} can be
obtained as a particular case of the one we have described
(essentially, it corresponds to a specific type of proximity graph
$\gG(n_1,n_2,\ldots,n_r)^+$). To apply his algorithm to a system of
multiplicities $\bold{m}=(m_0,m_1,\ldots,m_n)$, he uses successive
specializations, starting from a configuration of $n+1$ general
points of the plane and following with configurations corresponding
to the sequence of proximity graphs $\gG^1, \ldots, \gG^n$ where,
for each $k=1,2,\ldots,n$, $\gG^k$ is the one shown in Figure 3
(assuming that $\gG^1$ has no curved-dotted edge). The knowledge of
the dimensions of all complete linear systems on the surfaces
obtained by blowing-up at the points of whichever configuration
whose associated proximity graph is $\gG^n$ allows him to deduce,
using a similar reasoning to the one explained in the algorithm we
present here (but adapted to the above mentioned specific sequence
of specializations), an upper bound of the regularity of $\bold{m}$.
From this explanation, it is easy to deduce that applying Roé's
algorithm to $\bold{m}$ is equivalent to computing our bound
$\beta(\bold{m})$ taking the graph $\gG(n+1)^+$, adding previously
to $\bold{m}$ the suitable number of zeros. This proximity graph
corresponds, for instance, to the curve whose equation in projective
coordinates $(X:Y:Z)$ is $X Z^n+Y^{n+1}=0$, $Z=0$ being the line of
infinity.

\begin{figure}[h] \label{figure2}
\begin{center}
\setlength{\unitlength}{1mm}
\begin{picture}(25,60)

\qbezier[30](15,5)(5,17.5)(15,30)

\put(15,5){\circle*{1}} \put(15,5){\line(0,1){5}}
\put(15,10){\circle*{1}} \put(15,10){\line(0,1){5}}
\put(15,15){\circle*{1}} \put(15,15){\line(0,1){2.5}}

\put(15,30){\circle*{1}} \put(15,30){\line(0,1){2.5}}
\put(15,30){\line(0,-1){2.5}}

 \put(15,45){\line(0,-1){2.5}}

\put(15,45){\circle*{1}} \put(15,45){\line(0,1){5}}
\put(15,50){\circle*{1}} \put(15,50){\line(0,1){5}}
\put(15,55){\circle*{1}}

\put(14.7,21){$\vdots$}


\put(14.7,36){$\vdots$}

\put(16,5){ \scriptsize 0}\put(16,10){ \scriptsize 1}\put(16,15){
\scriptsize 2}\put(16,45){ \scriptsize $n-2$}\put(16,50){
\scriptsize $n-1$}\put(17,55){\scriptsize $n$}
\put(17,30){\scriptsize $k$}

\put(15,-2){$\gG^k$}


\end{picture}
\end{center}
\caption{Proximity graphs used in \cite{roe3}}
\end{figure}
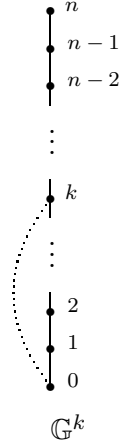

Note that, for each integer $n\geq 2$, the proximity graph
$\gG(\lfloor n/2 \rfloor+1)^+$ (considered in Example \ref{ex1}) is
the graph $\gG^k$ of Figure 3 for $k=\lfloor n/2 \rfloor$, if $n$ is
odd, and the one obtained from the same graph adding a new vertex
corresponding to a free point at the top, if $n$ is even. The fact
that this is one of the intermediate proximity graphs which appear
in the sequence of specializations used in \cite{roe3} and easy
reasonings concerning semicontinuity imply that our bounds
$\beta(\bold{m})$, taking the above proximity graph, are either
equal or lower than those obtained from \cite{roe3}. For homogeneous
systems of multiplicities $\bold{m}=(m,m,\ldots,m)$ and for a fixed
value of $n$, examples show that the difference between both bounds
increases when the multiplicity does so. For instance, this is the
behavior for $n+1=1000$ and $m$ taking values between 1 and 100. In
fact, when $m\leq 38$ the two bounds coincide, $\beta(\bold{m})$ is
sometimes better when $39\leq m\leq 68$ (in which case, the
difference is 1) and it is always better when $69\leq m\leq 100$
(the difference is 1 in all cases except for $m=98$, where it equals
2). Also, for $m=500$ (resp., $m=800$) (resp., $m=1200$),
$\beta(\bold{m})=16014$ (resp., $\beta(\bold{m})=25617$) (resp.,
$\beta(\bold{m})=38417$) and Roé's bound is $16021$ (resp., 25629)
(resp., $38436$). However, \cite{hr} gives better values in all the
checked cases where our bound is less than Roé's one.

For quasihomogeneous systems of multiplicities, there are cases in
which our bound seems to improve the existing ones (as far as the
author knows). As an example, consider the system of multiplicities
$\bold{m}=(4000,1000_{19})$ (where the subindex is the number of
occurrences). Taking the graph $\gG(10)$, it is obtained the bound
$\beta(\bold{m})=6009$. The Harbourne-Hirschowitz Conjecture
predicts that the regularity of $\bold{m}$ is 5917 and the bounds
provided in \cite{hir}, \cite{gim}, \cite{cat}, \cite{roe3},
\cite{bal} and \cite{xu} are 8367, 8000, 8000, 6183, 11140 and 6238,
respectively. Also, the bound 6015 is obtained by an algorithm based
on the reduction method described in \cite{dum1} and Cremona
transformations (it has been computed by using the computer program
provided in \cite{dumweb}), and the bound 7667 is obtained by using
the algorithm given in \cite{harbc} for computing the dimension of
line bundles on an smooth rational surface $X$ with anticanonical
bundle having an irreducible and reduced global section $D$, with
the further assumption that the morphism $\pic(X)\rightarrow
\pic(D)$ induced by the inclusion $D\subseteq X$ has trivial kernel.
The bound provided by the algorithm of \cite{hr}, using the
parameters $r=9$ and $d=2$, is 7667 (this algorithmic bound depends
on the choice of two parameters, but it is not clear how to obtain
the optimal values).

Now, consider the family of systems of multiplicities
$\bold{m}(m):=(m,1000_{19})$ for $m\geq 1000$. By applying
Corollary \ref{super2} to the graph $\gG(10)$ (see also Example
\ref{ex1}) it can be deduced that $\tau(\bold{m}(m))=m+1000$ when
$m\geq 9000$. Computing the above mentioned bounds of the
regularity for the remaining values of $m$ it holds that, when
either $m\in \{1619, 1622, 1623\}$ or $1625\leq m\leq 7765$, the
value $\beta(\bold{m}(m))$ (taking the graph $\gG(10)$) is less
than all the non-parametric bounds given in \cite{hir},
\cite{gim}, \cite{cat}, \cite{roe3}, \cite{bal} and \cite{xu}.
When either $3935\leq m\leq 3939$, $3944\leq m\leq 4081$,
$4083\leq m\leq 4085$ or $m\in \{3941,3942,4087,4089,4090,4092\}$
it holds that the bound $\beta(\bold{m}(m))$ is also better than
the one provided in \cite{dum1} and \cite{dumweb}; moreover, in
these cases, we have not found any pair of parameters $(r,d)$ for
which the bound given in \cite{hr} improves $\beta(\bold{m}(m))$.
It is worth adding that, by looking at systems of multiplicities
of the type $(m,h_{19})$ with $h\in \{1100,1200,1300,1400,1500\}$,
we have observed an increasing tendency (when $h$ grows) on the
number of values of $m$ for which $\beta(m,h_{19})$ seems to be
the best bound.

Although, in order to establish comparisons, it is natural to look
at homogeneous and quasihomogeneous cases, our algorithm can be
applied to arbitrary systems of multiplicities. Finally we notice
that, when the system of multiplicities $\bold{m}$ is either
homogeneous or quasihomogeneous, examples suggest that the bound
$\beta(\bold{m})$ is better when the graph $\gG(\lfloor n/2
\rfloor+1)^+$ is taken (where $n+1$ is the length of $\bold{m}$).

\begin{rem}
{\rm The results proved in Section \ref{infinity} and the
explanations given in the current section suggest that the algorithm
provided in \cite{roe2} for giving a lower bound of the least degree
$d$ such that a linear system ${\mathcal L}_d(\bold{m})$ is not
empty can also be generalized. However, we have not found evidences
of any significant improvement of this generalization with respect
to the existing bounds. }
\end{rem}

\vspace{3mm}
\par
\begin{center}
\footnotesize Dept. de Matemàtiques (ESTCE), UJI, Campus Riu Sec. \\
\footnotesize 12071  Castelló. SPAIN. \\ \footnotesize
monserra@mat.uji.es
\end{center}
\end{document}